\documentclass{gtart}

\def\ifplaintex{\expandafter\ifx\csname documentclass\endcsname\relax}

\ifplaintex 
\else
\fi

\def\gtm{{\mathsurround=0pt\it $\cal G\mskip-2mu$eometry \&\ 
$\cal T\!\!$opology $\cal M\mskip-1mu$onographs}}    

\def\gtp{{\mathsurround=0pt\it $\cal G\mskip-2mu$eometry \&\ 
$\cal T\!\!$opology $\cal P\!$ublications}}  

\def\recd{{\small Received:\qua\receiveddate\ifx\reviseddate\relax
\else\qquad Revised:\qua\reviseddate\fi\par}} 

\def\volumenumber#1{\def\thevolumenumber{#1}}
\def\volumeyear#1{\def\thevolumeyear{#1}}
\def\volumename#1{\def\thevolumename{#1}}
\def\papernumber#1{\def\thepapernumber{#1}}
\def\pagenumbers#1#2{\def\startpage{#1}\def\finishpage{#2}}
\def\published#1{\def\publishdate{#1}}
\def\received#1{\def\receiveddate{#1}}
\def\revised#1{\def\reviseddate{#1}}
\def\accepted#1{\def\accepteddate{#1}}

\long\def\asciiabstract#1{\long\def\theasciiabstract{#1}}

\let\\\par
\let\thevolumenumber\relax\let\thepapernumber\relax
\let\thevolumeyear\relax\let\startpage\relax
\let\finishpage\relax\let\publishdate\relax\let\receiveddate\relax
\let\reviseddate\relax\let\accepteddate\relax\let\theasciititle\relax
\let\theasciiauthors\relax
\let\theasciiabstract\relax

\let\theerratum\relax\let\theasciiemail\relax
\let\theshortauthors\relax\let\theshorttitle\relax

\def\startpage{1}\def\finishpage{15}\def\thepapernumber{77}

\volumenumber{2}
\volumename{Proceedings of the Kirbyfest}
\volumeyear{1999}

\long\def\maketitlep{   

\count0=\startpage

\gtm\nl        
{\small Volume \thevolumenumber: \thevolumename\nl 
\ifx\theerratum\relax\else Erratum \erratumnumber\nl\fi
Pages \startpage--\finishpage\nl}

\vglue 0.1truein   

{\parskip=0pt\leftskip 0pt plus 1fil\def\\{\par\smallskip}{\ifplaintex\large
\else\Large\fi\bf\thetitle}\par\medskip}   
\vglue 0.05truein 

{\parskip=0pt\leftskip 0pt plus 1fil\def\\{\par}{\sc\theauthors}
\par\medskip}%
 
\vglue 0.03truein 

{\small\leftskip 25pt\rightskip 25pt{\bf Abstract}\stdspace\theabstract

{\bf AMS Classification}\stdspace\theprimaryclass
\ifx\thesecondaryclass\relax\else; \thesecondaryclass\fi\par
{\bf Keywords}\stdspace \thekeywords\par}\vglue 7pt

}   

\font\phead=cmsl9 scaled 950
\font=cmsl9 scaled 1050
\font\pnum=cmbx10 scaled 913
\font\lnum=cmbx10 
\font\pfoot=cmsl9 scaled 950
\font=cmsl9 scaled 1050
\ifplaintex
\headline{\vbox to 0pt{\vskip -4.5mm\line{\small\phead\ifnum
\count0=\startpage ISSN 1464-8997 (on line)
1464-8989 (printed) \hfill {\pnum\folio}\else\ifodd\count0\def\\{ }%
\ifx\theshorttitle\relax\thetitle\else\theshorttitle\fi\hfill{\pnum\folio}
\else\def\\{ and }{\pnum\folio}\hfill\ifx\theshortauthors\relax\theauthors
\else\theshortauthors\fi\fi\fi}\vss}}
\footline{\vbox to 0pt{\vglue 0mm\line{\small\pfoot\ifnum\count0=\startpage
Published \publishdate:\qua\copyright\ \gtp\hfill\else
\gtm, Volume \thevolumenumber\ (\thevolumeyear)\hfill\fi}\vss
}}
\else
\makeatletter
\def\@oddhead{{\small\lhead\ifnum\count0=\startpage ISSN 1464-8997 (on line)
1464-8989 (printed) \hfill {\lnum\number\count0}\else\ifodd\count0
\def\\{ }\ifx\theshorttitle\relax \thetitle \else\theshorttitle\fi\hfill
{\lnum\number\count0}\else\def\\{ and }{\lnum\number\count0}
\hfill\ifx\theshortauthors\relax 
\theauthors\else\theshortauthors\fi\fi\fi}}\def\@evenhead{@oddhead}
\def\@oddfoot{\small\lfoot\ifnum\count0=\startpage Published \publishdate:\qua\copyright\ \gtp\hfill\else
\gtm, Volume \thevolumenumber\ (\thevolumeyear)\hfill\fi}
\def\@evenfoot{@oddfoot}
\makeatother
\fi

\let\maketitlepage\maketitlep
\let\makeshorttitle\maketitlepage
\let\maketitle\maketitlepage

\newwrite\gtoutfile
\long\gdef\makeheadfile{  
{\def\\{, }\def\s{ }
\immediate\openout\gtoutfile head.xxx
\immediate\write\gtoutfile{To: math@arxiv.org}
\immediate\write\gtoutfile{Subject: put OR rep NNNNN:ppppp}
\immediate\write\gtoutfile{--text follows this line--}
\immediate\write\gtoutfile{Proxy-for: \ifx\theasciiauthors\relax
\theauthors\else\theasciiauthors\fi\s<\ifx\theasciiemail\relax\theemail\else\theasciiemail\fi>}
\immediate\write\gtoutfile{\noexpand\\}
\immediate\write\gtoutfile{Authors: \ifx\theasciiauthors\relax
\theauthors\else\theasciiauthors\fi}
{\def\\{ }\immediate\write\gtoutfile{Title: \ifx\theasciititle\relax
\thetitle\else\theasciititle\fi}}
\immediate\write\gtoutfile{Subj-class: GT or SG, GR etc}
\immediate\write\gtoutfile{MSC-class: \theprimaryclass\ifx\thesecondaryclass\relax\else, \thesecondaryclass\fi}
\immediate\write\gtoutfile{Journal-ref: Geom. Topol. Monogr. \thevolumenumber\s
(\thevolumeyear) \startpage-\finishpage}
\immediate\write\gtoutfile{Comments: Published by Geometry and Topology Monographs at}
\immediate\write\gtoutfile{\s\s\s  http://www.maths.warwick.ac.uk/gt/GTMon\thevolumenumber/paper\thepapernumber.abs.html}
\immediate\write\gtoutfile{\noexpand\\}
\immediate\write\gtoutfile{}
\ifx\theasciiabstract\relax
\immediate\write\gtoutfile{\theabstract}\else
\immediate\write\gtoutfile{\theasciiabstract}\fi
\immediate\write\gtoutfile{}
\immediate\write\gtoutfile{\noexpand\\}
\immediate\write\gtoutfile{}
\immediate\closeout\gtoutfile}}  

\def\maketitlepage{\maketitlep\makeheadfile}
\let\makeshorttitle\maketitlepage
\let\maketitle\maketitlepage

\volumenumber{4}
\volumename{Invariants of knots and 3-manifolds (Kyoto 2001)}
\volumeyear{2002}
\papernumber{7}
\pagenumbers{89}{101}
\received{29 November 2001}
\revised{7 March 2002}
\accepted{22 July 2002}
\published{28 July 2002}

\usepackage{amsmath,amssymb}

\usepackage{eepic,epic}

\newtheorem{thm}{Theorem}[section]
\newtheorem{prop}[thm]{Proposition}
\newtheorem{coro}[thm]{Corollary}
\newtheorem{lem}[thm]{Lemma}

\theoremstyle{definition}
\newtheorem{defi}[thm]{Definition}
\newtheorem{que}[thm]{Question}
\newtheorem{exa}[thm]{Example}

\renewcommand{\theequation}
       {{\thesection}.\arabic{equation}}

\begin{document}
\title{Polynomial invariants and Vassiliev invariants}
\author{Myeong-Ju Jeong\\Chan-Young Park}

\address{Department of Mathematics, College of Natural
Sciences\\Kyungpook National University, Taegu 702-701 Korea}
\email{determiner@hanmail.net, chnypark@knu.ac.kr}

\begin{abstract}
We give a criterion to detect whether the derivatives of the
HOMFLY polynomial at a point is a Vassiliev invariant or not. In
particular, for a complex number $b$ we show that the derivative
$P_K^{(m,n)}(b,0)={{\partial^m}\over{\partial
a^m}}{{\partial^n}\over{\partial x^n}}P_K(a,x)|_{(a, x) = (b, 0)}$
of the HOMFLY polynomial of a knot $K$ at $(b,0)$ is a Vassiliev
invariant if and only if $b = \pm 1$. Also we analyze the space
$V_n$ of Vassiliev invariants of degree $\leq n$ for $n = 1, 2, 3,
4, 5$ by using the $~\bar{}~$--operation and the $^*$--operation
in \cite{JP2}. These two operations are unified to the
$~\hat{}~$--operation. For each Vassiliev invariant $v$ of degree
$\leq n$, $\hat{v}$ is a Vassiliev invariant of degree $\leq n$
and the value $\hat{v}(K)$ of a knot $K$ is a polynomial with
multi--variables of degree $\leq n$ and we give some questions on
polynomial invariants and the Vassiliev invariants.
\end{abstract}

\asciiabstract{We give a criterion to detect whether the derivatives
of the HOMFLY polynomial at a point is a Vassiliev invariant or
not. In particular, for a complex number b we show that the derivative
P_K^{(m,n)}(b,0)=d^m/da^m d^n/dx^n P_K(a,x)|(a, x) = (b, 0) of the
HOMFLY polynomial of a knot K at (b,0) is a Vassiliev invariant if and
only if b= -+1.  Also we analyze the space V_n of Vassiliev invariants
of degree <=n for n = 1,2,3,4,5 by using the bar-operation and the
star-operation in [M-J Jeong, C-Y Park, Vassiliev invariants and knot
polynomials, to appear in Topology and Its Applications].  These two
operations are unified to the hat-operation.  For each Vassiliev
invariant v of degree <=n, hat(v) is a Vassiliev invariant of degree
<=n and the value hat(v)K) of a knot K is a polynomial with
multi-variables of degree <=n and we give some questions on polynomial
invariants and the Vassiliev invariants.}



\maketitlepage

\let\\\par

\section{Introduction}
In 1990, V. A. Vassiliev introduced the concept of a finite type
invariant of knots, called Vassiliev invariants \cite{V}. There
are some analogies between Vassiliev invariants and polynomials.
For example, in 1996 D. Bar--Natan showed that when a Vassiliev
invariant of degree $m$ is evaluated on a knot diagram having $n$
crossings, the result is approximately bounded by a constant times
of $n^m$ \cite{BN2} and S. Willerton \cite{W2} showed that for
any Vassiliev invariant $v$ of degree $n$, the function
$p_v(i,j)\co =v(T_{i,j})$ is a polynomial of degree $\leq n$ for
each variable $i$ and $j$. Recently, we \cite{JP1} defined a
sequence of knots or links induced from a double dating tangle and
showed that any Vassiliev invariant has a polynomial growth on
this sequence.

J. S. Birman and X.--S. Lin \cite{BL} showed that each
coefficient in the Maclaurin series of the Jones, Kauffman, and
HOMFLY polynomial, after a suitable change of variables, is a
Vassiliev invariant, and T. Kanenobu \cite{Ka2,Ka3}
showed that some derivatives of the HOMFLY and the Kauffman
polynomial are Vassiliev invariants. For the question whether the
$n$--th derivatives of knot polynomials are Vassiliev invariants
or not, we \cite{JP2} gave complete solutions for the Jones,
Alexander, Conway polynomial and a partial solution for the
$Q$--polynomial. Also we introduced the $~\bar{}~$--operation and
the $^*$--operation to obtain polynomial invariants from a
Vassiliev invariant of degree $n$. From each of these new
polynomial invariants, we may get at most $(n+1)$ linearly
independent numerical Vassiliev invariants.

In this paper, we find a line and two points in the complex plane
where the derivatives of the HOMFLY polynomial can possibly be
Vassiliev invariants and analyze the space $V_n$ of Vassiliev
invariants for $n \leq 5$ by using the $~\bar{}~$--operation and
the $^*$--operation.

Throughout this paper all knots or links are assumed to be
oriented unless otherwise stated. For a knot $K$ and $i\in
\mathbb{N}$, $K^i$ denotes the $i$--times self--connected sum of
$K$ and $\mathbb{N},$ $\mathbb{Z},$ $\mathbb{Q},$ $\mathbb{R},$
$\mathbb{C}$ denote the sets of nonnegative integers, integers,
rational numbers, real numbers and complex numbers, respectively.

A knot or link invariant $v$ taking values in an abelian group can
be extended to a singular knot or link invariant by taking the
difference between the positive and negative resolutions of the
singularity. A knot or link invariant $v$ is called a {\it
Vassiliev invariant of degree $n$} if $n$ is the smallest
nonnegative integer such that $v$ vanishes on singular knots or
links with more than $n$ double points. A knot or link invariant
$v$ is called a {\it Vassiliev invariant} if $v$ is a Vassiliev
invariant of degree $n$ for some nonnegative integer $n$.

\begin{defi}\cite {JP1}\qua
Let {\bf J} be a closed interval $[a,b]$ and $k$ a positive
integer. Fix $k$ points in the upper plane {\bf J}$^2\times \{b\}$
of the cube {\bf J}$^3$ and their corresponding $k$ points in the
lower plane {\bf J}$^2\times \{a\}$ of the cube {\bf J}$^3$. A
$(k,k)$--{\it tangle} is obtained by attaching, within {\bf
J}$^3$, to these $2k$ points $k$ curves, none of which should
intersect each other. A $(k, k)$--tangle is said to be {\it
oriented} if each of its $k$ curves is oriented. Given two
$(k,k)$--tangles $S$ and $T$, roughly the {\it tangle product}
$ST$ is defined to be the tangle obtained by gluing the lower
plane of the cube containing $S$ to the upper plane of the cube
containing $T$. The {\it closure} $\overline T$ of a tangle $T$ is
the unoriented knot or link obtained by attaching $k$ parallel
strands connecting the $k$ points and their corresponding $k$
points in the exterior of the cube containing T. When the tangles
$S$ and $T$ are oriented, the oriented tangle $ST$ is defined only
when it respects the orientations of $S$ and $T$ and the closure
$\overline S$ has the orientation inherited from that of $S$ and
$\overline {ST}$ is the oriented knot or link obtained by closing
the $(k, k)$--tangle $ST$.
\end{defi}

\begin{defi}\cite {JP1}\qua
An oriented $(k,k)$--tangle $T$ is called a {\it double dating
tangle} ({\it DD--tangle} for short) if there exist some ordered
pairs of crossings of the form $(*)$ in Figure \ref{bistar}, so
that $T$ becomes the trivial $(k,k)$--tangle when we change all
the crossings in the ordered pairs, where $i$ and $j$ in Figure
\ref{bistar}, denote components of the tangle. Note that a
DD--tangle is always an oriented tangle.
\end{defi}

\begin{figure}[ht!]
\begin{center}
\setlength{\unitlength}{0.00041667in}
\begingroup\makeatletter\ifx\SetFigFont\undefined%
\gdef\SetFigFont#1#2#3#4#5{%
  \reset@font\fontsize{#1}{#2pt}%
  \fontfamily{#3}\fontseries{#4}\fontshape{#5}%
  \selectfont}%
\fi\endgroup%
{\renewcommand{\dashlinestretch}{30}
\begin{picture}(3618,2484)(0,-10)
\put(1509.000,1350.000){\arc{3000.000}{2.4981}{3.7851}}
\put(2109.000,1350.000){\arc{3000.000}{5.6397}{6.9267}}
\path(609,450)(1509,2250)
\path(1509,450)(1134,1200)
\path(609,2250)(984,1500)
\path(2109,450)(3009,2250)
\path(3009,450)(2634,1200)
\path(2109,2250)(2484,1500)
\path(609,2100)(609,2250)(759,2175)
\path(1359,2175)(1509,2250)(1509,2100)
\path(2859,2175)(3009,2250)(3009,2100)
\path(2859,525)(3009,450)(3009,600)
\put(534,2325){\makebox(0,0)[lb]{\smash{{{\SetFigFont{8}{7.2}{\rmdefault}{\mddefault}{\updefault}$i$}}}}}
\put(1434,2325){\makebox(0,0)[lb]{\smash{{{\SetFigFont{8}{7.2}{\rmdefault}{\mddefault}{\updefault}$j$}}}}}
\put(2034,2325){\makebox(0,0)[lb]{\smash{{{\SetFigFont{8}{7.2}{\rmdefault}{\mddefault}{\updefault}$i$}}}}}
\put(2934,2325){\makebox(0,0)[lb]{\smash{{{\SetFigFont{8}{7.2}{\rmdefault}{\mddefault}{\updefault}$j$}}}}}
\put(1809,450){\makebox(0,0)[lb]{\smash{{{\SetFigFont{8}{7.2}{\rmdefault}{\mddefault}{\updefault}$,$}}}}}
\end{picture}
}
\end{center}
\vglue -4mm
\caption{$(*)$}\label{bistar}
\end{figure}

Since every $(1,1)$--tangle is a double dating tangle, every knot
is a closure of a double dating $(1,1)$--tangle. But there is a
link which is not the closure of any DD--tangle since the linking
number of two components of the closure of a DD--tangle must be
$0.$

\begin{defi}\cite {JP1}\qua
Given an oriented $(k,k)$--tangle $S$ and a double dating
$(k,k)$--tangle $T$ such that the product $ST$ is well--defined,
we have a sequence of links $\{L_i(S,T)\}_{i=0}^{\infty}$ obtained
by setting $L_i(S,T)=\overline {ST^i}$ where $T^i=TT\cdots T$ is
the $i$--times self--product of $T$ and $T^0$ is the trivial $(k,
k)$--tangle. We call $\{L_i(S,T)\}_{i=0}^{\infty}$
($\{L_i\}_{i=0}^{\infty}$ for short) the {\it sequence induced
from the $(k, k)$--tangle $S$ and the double dating $(k,
k)$--tangle $T$} or simply a {\it sequence induced from the double
dating tangle T}.
\end{defi}

In particular, if $\overline{S}$ is a knot for a $(k, k)$--tangle
$S$, then $L_i(S,T)=\overline {ST^i}$ is a knot for each $i \in
\mathbb{N}$ since $T^i$ can be trivialized by changing some
crossings.

\begin{thm}\label{vasddt}{\rm\cite{JP2}}\qua
Let $\{L_i\}_{i=0}^\infty$ be a sequence of knots induced from a
DD--tangle. Then any Vassiliev knot invariant $v$ of degree $n$
has a polynomial growth on $\{L_i\}_{i=0}^\infty$ of degree $\leq
n.$
\end{thm}

\begin{coro}\label{vascon}{\rm\cite{JP2}}\qua
Let $L$ and $K$ be two knots. For each $i \in \mathbb{N}$, let
$K_i=K\sharp L\sharp\cdots\sharp L$ be the connected sum of $K$ to
the $i$--times self--connected sum of $L$. If $v$ is a Vassiliev
invariant of degree $n$, then $v|_{\{K_i\}_{i=0}^{\infty}}$ is a
polynomial function in $i$ of degree $\leq n$.
\end{coro}

The converse of Corollary \ref{vascon} is not true. In fact, the
maximal degree $u(K)$ of the Conway polynomial $\nabla_K(z)$ for a
knot $K$ is a counterexample.

\rk{Acknowledgement} This work was supported by grant No.\
R02-2000-00008 from the Basic Research Program of the Korea
Science Engineering Foundation.

\section{The derivatives of the HOMFLY polynomial and Vassiliev
invariants.} From now on, the notations $3_1$, $4_1$, $5_1$ and
$6_1$ will mean the knots in the Rolfsen's knot table \cite{R}.
For the definitions of the HOMFLY polynomial $P_L(a,z)$ and the
Kauffman polynomial $F_L(a,x)$ of a knot or link $L$, see
\cite{Kaw}.

Note that the {\it Jones polynomial} $J_L(t)$, the {\it Conway
polynomial} $\nabla_L(z)$, and the {\it Alexander polynomial}
$\Delta_L(t)$ of a knot or link $L$ can be defined from the HOMFLY
polynomial $P_L(a,z) \in {\mathbb Z}[a,a^{-1},z,z^{-1}]$ via the
equations $J_L(t)= P_L(t,t^{1/2}-t^{-1/2})$,
$\nabla_L(z)=P_L(1,z)$ and $\Delta_L(t)=P_L(1,t^{1/2}-t^{-1/2})$
respectively and that the {\it $Q$--polynomial} $Q_L(x)$ can be
defined from the Kauffman polynomial $F_L(a,x)$ via the equation
$Q_L(x) = F_L(1,x)$.

By using the skein relations, we can see that $P_L(a,z)$ and
$F_L(a,x)$ are {\it multiplicative under the connected sum.} i.e.\
$P_{L_1\sharp L_2}(a,z) = P_{L_1}(a,z)P_{L_2}(a,z)$ and
$F_{L_1\sharp L_2}(a,x) = F_{L_1}(a,x)F_{L_2}(a,x)$ for all knots
or links $L_1$ and $L_2$. So the Jones, Conway, Alexander and
$Q$--polynomials are also multiplicative under the connected sum.

It is well known that $P_K(a,z)\in {\mathbb Z}[a^2,a^{-2},z^2]$
and  $F_K(a,x) \in {\mathbb Z}[a,a^{-1},x]$ for a knot $K$. For
each $i\in \mathbb{N}$ and each knot $K$, we denote by $F_i(K; a)$
and $P_{2i}(K; a)$ the coefficient of $x^i$ in $F_K(a, x)$ and the
coefficient of $z^{2i}$ in $P_K(a,z),$ respectively, which are
polynomials in $a$.

Throughout this section, knot polynomials are always assumed to be
multiplicative under the connected sum.

We consider 1--variable knot polynomials first and then
2--variable knot polynomials.

\begin{lem}\label{derpol1}{\rm\cite{JP2}}\qua
Let $f_K(x)$ be a knot polynomial of a knot $K$ such that $f_K(x)$
is infinitely differentiable in a neighborhood of a point $a$ and
assume that $f_K^{(1)}(a)\neq 0$. Then there exists a unique
polynomial $p(x)$ of degree $m$ such that
$f_{K^i}^{(m)}(a)=(f_K(a))^ip(i)$ for $i>m$.
\end{lem}
\eject
\begin{thm}\label{polvas1} {\rm\cite{JP2}}\qua
 For each $n \in \mathbb{N}$, we have\\
$\mathrm{(1)}$\qua $J_K^{(n)}(a)$ is a Vassiliev invariant if and only
if $a=1$.
\\ $\mathrm{(2)}$\qua $\nabla_K^{(n)}(a)$ is a Vassiliev invariant if and only if
$a=0$.
\\ $\mathrm{(3)}$\qua $\Delta_K^{(n)}(a)$ is a Vassiliev invariant if and only if
$a=1$. \\ $\mathrm{(4)}$\qua $Q_K^{(n)}(a)$ is not a Vassiliev
invariant if $a\neq -2, 1$.
\end{thm}

\begin{thm}\label{change}
Let $g\co {\mathbb R} \rightarrow {\mathbb R}$ be infinitely
differentiable function at $x=a$ with $g^{(1)}(a)\neq 0$. Assume
that $f_K(x)$ is a knot polynomial which is infinitely
differentiable in a neighborhood of $g(a)$ for all knots $K$ and
that there exists a knot $L$ such that $f_L(g(a))\neq 0,1$ and
$f_L^{(1)}(g(a))\neq 0$. Then each coefficient of $(x-a)^n$ in the
Taylor expansion of $f_K\circ g(x)$ at $x = a,$ is not a Vassiliev
invariant. \end{thm}
\begin{proof} Consider a sequence $\{L^i\}_{i=0}^{\infty}$ of knots.
By Lemma \ref{derpol1}, we see that $(f_{L^i}(g(x)))^{(n)}|_{x=a}
= (f_L(g(a)))^ip(i),$ where $p(i)$ is a polynomial in $i$ of
degree $n,$ and hence the coefficient
${{1}\over{n!}}(f_{K}(g(x)))^{(n)}|_{x=a}$ of $(x-a)^n$ does not
have a polynomial growth on $\{L^i\}_{i=0}^{\infty}$.

It follows from Corollary \ref{vascon} that the coefficient of
$(x-a)^n$ in the Taylor expansion of $f_K\circ g(x)$ is not a
Vassiliev invariant.
\end{proof}

J. S. Birman and X.--S. Lin \cite{BL} showed that each
coefficient in the Maclaurin series of $J_K(e^x)$ is a Vassiliev
invariant. As a generalization of Birman and Lin's type of
changing variables, we have

\begin{thm}\label{coeff} Let $g\co {\mathbb R} \rightarrow
{\mathbb R}$ be an infinitely differentiable
function at $x=a.$ Assume that $g^{(1)}(a)\neq 0$. Then \\
$\mathrm{(1)}$\qua each coefficient of $(x-a)^n$ in the Taylor
expansion of $J_K\circ g(x)$ at $x = a,$ is a Vassiliev invariant
if and only if $g(a) = 1,$\\
$\mathrm{(2)}$\qua each coefficient of $(x-a)^n$ in the Taylor
expansion of $\nabla_K\circ g(x)$ at $x = a,$ is
a Vassiliev invariant if and only if $g(a) = 0,$\\
$\mathrm{(3)}$\qua each coefficient of $(x-a)^n$ in the Taylor
expansion of $\Delta_K\circ g(x)$ at $x = a,$ is
not a Vassiliev invariant if and only if $g(a) = 1$ and \\
$\mathrm{(4)}$\qua if $g(a)\neq -2, 1$ then each coefficient of
$(x-a)^n$ in the Taylor expansion of $Q_K\circ g(x)$ at $x = a,$
is not a Vassiliev invariant.
\end{thm}
\begin{proof}
(1)\qua Let $A_K = \{t|~J_K(t)=0,1\}\bigcup\{t|~J_K^{(1)}(t)=0\}$ for
a knot $K$. Then $A_{3_1}\bigcap A_{4_1} = \{ 1\}$. Thus if
$g(a)\neq 1$, then $g(a)\in \mathbb{R} \setminus (A_{3_1}\bigcap
A_{4_1})$. Take $L = 3_1$ in Theorem \ref{change} if $g(a)\in
\mathbb{R} \setminus A_{3_1}$ and $L = 4_1$ in Theorem
\ref{change} if $g(a)\in \mathbb{R} \setminus A_{4_1}.$ Then
$J_L(g(a))\neq 0,1$ and $J_L^{(1)}(g(a))\neq 0$. So by Theorem
\ref{change}, each coefficient of $(x-a)^n$ in the Taylor
expansion of $J_K\circ g(x)$ is not a Vassiliev invariant.
Conversely, assume that $g(a)=1$ and that $n \in \mathbb{N}$.
Since the coefficient of $(x-a)^n$ in the Taylor expansion of
$J_K(g(x))$ is a linear combination of $1, J_K^{(1)}(1), \cdots,
J_K^{(n)}(1)$, by Theorem \ref{polvas1}, it is a Vassiliev
invariant. The proofs of (2), (3) and (4) are similar.
\end{proof}

\begin{exa} Take $f(x) = \mathrm{sin}(x)$ for $x \in
\mathbb{R}.$ Then $f(0)\neq 1$ and $f^{(1)}(0)\neq 0$. Thus each
coefficient in the Maclaurin series of $J_K(\mathrm{sin}(x)) =
J_K(f(x))$ is not a Vassiliev invariant. But each coefficient in
the Maclaurin series of $\nabla_K(\mathrm{sin}(x)) =
\nabla_K(f(x))$ is a Vassiliev invariant, since it is a finite
linear combination of the coefficients of the Conway polynomial
$\nabla_K(z)$ of a knot $K$.
\end{exa}

\indent Now we will deal with 2--variable knot polynomials such as
the HOMFLY polynomial $P_K(a,z)\in {\mathbb Z}[a,a^{-1},z]$ and
the Kauffman polynomial $F_K(a,x)\in {\mathbb Z}[a,a^{-1},x]$. For
a 2--variable Laurent polynomial $g(x,y)$ which is infinitely
differentiable on a neighborhood of $(a,b)$, we denote
${{\partial^m}\over{\partial x^m}}{{\partial^n}\over{\partial
y^n}}g(x,y)|_{(x, y) = (a, b)}$ by $g^{(m,n)}(a,b)$ for each pair
$(m,n) \in \mathbb{N}^2$.

\begin{thm}\label{vasnes2}{\rm\cite{JP2}}\qua
Let $g_K(x,y)$ be a 2--variable knot polynomial which is
infinitely differentiable on a neighborhood of $(a,b)$ for all
knots $K$. If there exists a knot $L$ such that $g_L(a,b)\neq
0,1$, $g_L^{(1,0)}(a,b)\neq 0$  and $g_L^{(0,1)}(a,b)\neq 0$ then
$g_K^{(m,n)}(a,b)$ is not a Vassiliev invariant for all $m, n \in
\mathbb{N}$.
\end{thm}

\begin{lem} \label{derpol3}
 Let $g_K(x,y)$ be a 2--variable knot polynomial which is infinitely
differentiable on a neighborhood of $(a,b)\in \mathbb{C}^2$ for
all knots $K$ and let $m, n \in \mathbb{N}$. If there exists a
knot $L$ such that $g_L(a,b)\neq 0,1$, $g_L^{(1,0)}(a,b)\neq 0$,
$g_L^{(0,1)}(a,b)= 0$ and $g_L^{(0,2)}(a,b) \neq 0$ then there
exists a polynomial $p(i)$ of degree $m+n$ such that
$g_{L^i}^{(m,2n)}(a,b) = (g_L(a,b))^ip(i)$ for $i>m+2n$.
\end{lem}

\begin{proof} It is similar to that of Theorem 2.12 in \cite{JP2}.
\end{proof}

\begin{lem}\label{vascon3}
Let $g_K(x,y)$ be a 2--variable knot polynomial which is
infinitely differentiable on a neighborhood of $(a,b)\in
\mathbb{C}^2$ for all knots $K$. If there exists a knot $L$ such
that $g_L(a,b)\neq 0,1$, $g_L^{(1,0)}(a,b)\neq 0$,
$g_L^{(0,1)}(a,b)= 0$ and $g_L^{(0,2)}(a,b) \neq 0$ then
$g_K^{(m,2n)}(a,b)$ is not a Vassiliev invariant for all $m, n \in
\mathbb{N}$.
\end{lem}

\begin{proof} It follows from Lemma \ref{derpol3} and Corollary
\ref{vascon}.
\end{proof}

\begin{thm}\label{HOM}
Let $n\in \mathbb{N}$ and $a\in \mathbb{C}$. $P_{2i}^{(n)}(K;a)$
is a Vassiliev invariant if and only if $a = \pm 1$.
\end{thm}

\begin{proof} Note that $P_{2i}^{(n)}(K;a) =
(2i)!P_K^{(n,2i)}(a,0).$ Since $P_K(a,z) \in
\mathbb{Z}[a^2,a^{-2},z^2]$ for all knots $K$,
$P_K^{(n,1)}(a,0)=0$ for all $a\in \mathbb{C}$ and all knots $K$.
For each knot $K$, let $A_K^1 = \{a \in \mathbb{C}
~|~P_K(a,0)=0~{\text{or}}~1\}$, $A_K^2 = \{a \in \mathbb{C}
~|~P_K^{(1,0)}(a,0)=0\}$, $A_K^3 = \{a \in \mathbb{C}
~|~P_K^{(0,2)}(a,0)=0\}$ and $A_K=A_K^1\bigcup A_K^2 \bigcup
A_K^3.$ Since $P_{3_1}(a,z) = (-a^{-4}+2a^{-2})+a^{-2}z^2$ and
$P_{4_1}(a,z) = (a^{-2}-1+a^2)-z^2$, we have $A_{3_1}=\{\pm
{{\sqrt{2}}\over{2}}, \pm 1 \},$ $A_{4_1} =
\{\pm({{\sqrt{3}+\sqrt{-1}}\over{2}}),
\pm({{\sqrt{3}-\sqrt{-1}}\over{2}}), \pm 1, \pm \sqrt{-1}\}$ and
hence $A_{3_1}\bigcap A_{4_1} = \{\pm 1\}.$ Thus if $a\neq \pm 1$,
then, by Lemma \ref{vascon3}, $P_{2i}^{(n)}(K;a)$ is not a
Vassiliev invariant. Conversely, T. Kanenobu \cite{Ka3} showed
that $P_{2i}^{(n)}(K;1)$ is a Vassiliev invariant. Since
$P_{2i}^{(n)}(K;-1) = (-1)^nP_{2i}^{(n)}(K;1)$,
$P_{2i}^{(n)}(K;-1)$ is also a Vassiliev invariant.
\end{proof}

By Theorem \ref{HOM}, for $b\in \mathbb{C}$, $P_K^{(m,n)}(b,0)$ is
a Vassiliev invariant if and only if $n$ is odd or $b = \pm 1.$
For $(b,y)\in \mathbb{C}^2$ with $y\neq 0$, we have the following

\begin{thm}
Let $m, n$ be nonnegative integers. If $(b, y) \in \mathbb{C}^2$
with $y \neq 0$ such that $P_K^{(m,n)}(b,y)$ is a Vassiliev
invariant, then $(b,y) = (b, \pm(b-b^{-1})),$ $(\pm \sqrt{-1},
\sqrt{-3})$ or $(\pm \sqrt{-1}, -\sqrt{-3}).$
\end{thm}

\begin{proof} By direct calculations, $P_{3_1}(a,z) =
(-a^{-4}+2a^{-2})+a^{-2}z^2$, $P_{4_1}(a,z) = (a^{-2}-1+a^2)-z^2$
and $P_{6_1}(a,z) = (a^{-4}-a^{-2}+a^2)+z^2(-a^{-2}-1).$ Let
$A_K^1 = \{(b,y)~|~P_K(b,y)=0~{\text{or}}~1\}$, $A_K^2 =
\{(b,y)~|~P_K^{(1,0)}(b,y)=0\}$,  $A_K^3 =
\{(b,y)~|~P_K^{(0,1)}(b,y)=0\}$ and $A_k = A_K^1\bigcup A_K^2
\bigcup A_K^3$ for each knot $K$. Then
\begin{eqnarray*}
A_{3_1}\cap A_{4_1} &=& (A_{3_1}^1\cap A_{4_1}^1)\cup
(A_{3_1}^1\cap A_{4_1}^2)\cup \cdots \cup(A_{3_1}^3
\cap A_{4_1}^3)\\
&=& \{~(\pm \sqrt{-1}, 2\sqrt{-1}),~(\pm \sqrt{-1},
-2\sqrt{-1})\}\\
& &\cup \{(\pm \sqrt{-1}, \sqrt{-3}),~(\pm \sqrt{-1},
-\sqrt{-3}),~(\pm 1,\sqrt{-1}),~(\pm 1,-\sqrt{-1})\}\\
& &\cup\{({{-1 \pm \sqrt{5}}\over {2}},\sqrt{1 \pm
\sqrt{5}}),~({{-1 \pm \sqrt{5}}\over {2}},-\sqrt{1 \pm
\sqrt{5}})\}\\
& &\cup \{(b,y)~|~y = \pm(b-b^{-1})\}.
\end{eqnarray*}
So we get
\begin{eqnarray*}
&&A_{3_1}\cap A_{4_1}\cap A_{6_1}\\ &=& ((A_{3_1}\cap A_{4_1})\cap
A_{6_1}^1)\cup ((A_{3_1}\cap A_{4_1})\cap
A_{6_1}^2)\cup ((A_{3_1}\cap A_{4_1})\cap A_{6_1}^3)\\
&=& \{(b,y)~|~y = \pm(b-b^{-1})\} \cup \{(\pm \sqrt{-1},
\sqrt{-3}),~(\pm \sqrt{-1}, -\sqrt{-3})\}.
\end{eqnarray*}

If $(b,y)\in \mathbb{C}^2 \setminus (A_{3_1}\cap A_{4_1}\cap
A_{6_1})$, then, by Theorem \ref{vasnes2}, $P_K^{(m,n)}(b,y)$ is
not a Vassiliev invariant.
\end{proof}

Whether a finite product of the derivatives of knot polynomials at
some points is a Vassiliev invariant or not can be detected by
using Lemma \ref{derpol1}, Theorem \ref{vasnes2}, Lemma
\ref{derpol3} and Corollary \ref{vascon}. For example if there is
a knot $L$ such that $J_L^{(1)}(a) \neq 0, Q _L^{(1)}(b)\neq 0,
P_L^{(1,0)}(c,y) \neq 0, P_L^{(0,1)}(c,y) \neq 0$ and
$J_L(a)Q_L(b)P_L(c,y)$ $\neq 0, 1,$
 then the product
$J_K^{(k)}(a)Q_K^{(l)}(b)P_K^{(m,n)}(c,y)$ is not a Vassiliev
invariant for any $k, l, m, n \in \mathbb{N}$.

Since $Q_K^{(1)}(-2)=J_K^{(2)}(1)$ (T. Kanenobu \cite{Ka1}),
$Q_K^{(1)}(-2)$ is a Vassiliev invariant of degree $\leq 2.$ Note
that $Q_K^{(0)}(1) = 1$ for any knot $K$ and hence $Q_K^{(0)}(1)$
is a Vassiliev invariant of degree 0, but $Q_K^{(1)}(1)$ and
$Q_K^{(2)}(1)$ are not Vassiliev invariants \cite{JP2}.

\medskip

{\bf Open Problem} (A. Stoimenow
\cite{Sto1})\qua Is $Q_K^{(n)}(-2)$ a Vassiliev invariant for $n\geq
2$ ?

\begin{que} Is $Q_K^{(n)}(1)$ a Vassiliev invariant for $n \geq
3$ ?
\end{que}

The above two problems are the only remaining unsolved problems in
one variable knot polynomials \cite{JP2}.

\begin{que} Find all the points at which the
derivatives of the Kauffman polynomial are Vassiliev invariants.
\end{que}

\begin{que} Find all linear combinations of any finite
products of derivatives of knot polynomials, which are Vassiliev
invariants. \end{que}

\section{New polynomial invariants from Vassiliev invariants}
In this section, a Vassiliev invariant $v$ always means a
Vassiliev invariant taking values in a numerical number field {\bf
F} $=$ $\mathbb{Q},$ $\mathbb{R},$ or $\mathbb{C}$. We begin with
introducing the constructions of new polynomial invariants from a
given Vassiliev invariant (see \cite{JP1}) and then we will define
a new polynomial invariant unifying the polynomial invariants
obtained from the constructions in \cite{JP1}. The new polynomial
invariant is also a Vassiliev invariant and so we get various
numerical Vassiliev invariants from the coefficients of the new
polynomial invariant.

Let $K$ and $L$ be two knots and let $\{L_i\}_{i=0}^{\infty}$ be a
sequence of knots induced from a DD--tangle. Since any
$(1,1)$--tangle is a DD--tangle, we get two sequences $\{L\sharp
K^i\}_{i=0}^{\infty}$ and $\{K\sharp L_i\}_{i=0}^{\infty}$ of
knots induced from DD--tangles.

Let $v$ be a Vassiliev invariant of degree $n$ and fix a knot $L$.
Then by Corollary \ref{vascon}, for each knot $K$ there exist
unique polynomials $p_K(x)$ and $q_K(x)$ in ${\bf F}[x]$ with
degrees $\leq n$ such that $v(L\sharp K^i)=p_K(i)$ and $v(K\sharp
L_i)=q_K(i)$. We define two polynomial invariants $\bar v$ and
$v^*$ as follows: $\bar v\co \{\text{knots}\}\rightarrow {\bf
F}[x]$ by $\bar v(K)=p_K(x)$ and $v^*\co
\{\text{knots}\}\rightarrow {\bf F}[x]$ by $v^*(K)=q_K(x)$. Then
$\bar v(K)|_{x=j}=p_K(j)=v(L\sharp K^j)$ and
$v^*(K)|_{x=j}=q_K(j)=v(K\sharp L_j)$ for all $j\in \mathbb{N}$.

Then we have the following

\begin{thm}\label{vaspolth}{\rm\cite{JP2}}\qua
Let $v$ be a Vassiliev invariant of degree $n$ taking values in a
numerical field {\bf F}. \\ $\mathrm{(1)}$\qua For a fixed knot $L$,
$\bar v$ is a Vassiliev invariant of degree $\leq n$ and the
degree of $x$ in $\bar v(K)$ is $\leq n$. In particular if $L$ is
the unknot, $\bar v$ is a Vassiliev invariant of degree $n$ and
$\bar v(K)|_{x=1}=v(K)$.\\ $\mathrm{(2)}$\qua For a fixed sequence
$\{L_i\}_{i=0}^\infty$ of knots induced from a DD--tangle, $v^*$
is a Vassiliev invariant of degree $\leq n$ and the degree of $x$
in $v^*(K)$ is $\leq n$. In particular if $L_j$ is the unknot for
some $j \in \mathbb{N}$, then $v^*$ is a Vassiliev invariant of
degree $n$ and $v^*(K)|_{x=j}=v(K)$.
\end{thm}

Given a Vassiliev invariant $v$ of degree $n$, we may get at most
$(n+1)$ linearly independent numerical Vassiliev invariants which
are the coefficients of the polynomial invariants $\bar v$ and
$v^*$ respectively and then apply
 $~\bar{}~$--operation and $^*$--operation repeatedly on these new Vassiliev
invariants to get another new Vassiliev invariants. Inductively we
may obtain various Vassiliev invariants.

We note that for a Vassiliev invariant $v$ of degree $n$, since
$\bar v(K)$ and $v^*(K)$ are polynomials of degrees $\leq n$ for
any knot $K$, the polynomial invariants $\bar v$ and $v^*$ are
completely determined by $\{\bar v(K)|_{x=i}~|~0 \leq i \leq n\}$
and $\{v^*(K)|_{x=i}~|~0 \leq i \leq n\}$ respectively.

Let $V_n$ be the space of Vassiliev invariants of degrees $\leq n$
and let $A_n \subset V_n$. For each nonnegative integer $j$,
define $A_n^j$ as follows. Set $A_n^0=A_n$ and define inductively
$A_n^j$ to be the set of all Vassiliev invariants obtained from
the coefficients of the new polynomial invariants $\bar v$ and
$v^*$ ranging over all $v \in A_n^{j-1}$, all knots $L$ and all
sequences $\{L_i\}_{i=0}^{\infty}$ induced from all DD--tangles in
Theorem \ref{vaspolth}.

Define $A_n^* = \cup_{j=0}^{\infty}A_n^j$. We ask ourselves the
following:

\medskip

{\bf Question}\qua \cite{JP2}\qua Find
a minimal finite subset $A_n$ of
$V_n$ such that span($A_n^*)$ $= V_n$.\\

Let $V_n$ be the space of Vassiliev invariants of degree $\leq n$.
Then the dimension of $V_n/{V_{n-1}}$ is $0, 1, 1, 3, 4, 9, 14$
for $n = 1, 2, 3, 4, 5, 6, 7$ \cite{BN}.

\begin{prop}{\rm\cite{Ka2,Ka3}}\qua
For each nonnegative integer
$k$ and $l$,\\
$\mathrm{(1)}$\qua $P_{2k}^{(l)}(K;1)$ is a Vassiliev invariant of
degree $\leq
2k+l$.\\
$\mathrm{(2)}$\qua $(\sqrt{-1})^{k+l}F_k^{(l)}(K;\sqrt{-1})$ is a
Vassiliev invariant of degree $\leq k+l$.
\end{prop}

If $v_n$ and $v_m$ are Vassiliev invariants of degrees $n$ and $m$
respectively, then the product $v_nv_m$ is a Vassiliev invariant
of degree $\leq n+m$ \cite{BN,W}.

We get a base for each $V_n$ ($n\leq 5$) from the results of J. S.
Birman and X.--S. Lin (cite{BL}, D. Bar--Natan \cite{BN} and
T. Kanenobu \cite{Ka4}.
\begin{thm} \label{basis}
{\rm\cite{Ka4,BL,BN}}\qua Let $V_n$ be the
space of Vassiliev invariants of degree $\leq n$.
Then\\
$\mathrm{(1)}$\qua $\{1\}$ is a basis for $V_0 = V_1$, where $1$ is
the constant
map with image $\{1\}$.\\
$\mathrm{(2)}$\qua $\{a_2(K)\}$ is a basis for $V_2/{V_1}.$\\
$\mathrm{(3)}$\qua $\{J_K^{(3)}(1)\}$ is a basis for $V_3/{V_2}.$\\
$\mathrm{(4)}$\qua $\{(a_2(K))^2, a_4(K), J_K^{(4)}(1)\}$ is a basis for $V_4/{V_3}.$\\
$\mathrm{(5)}$\qua $\{a_2(K)P_0^{(3)}(K;1), P_0^{(5)}(K;1),
P_4^{(1)}(K;1), \sqrt{-1}F_4^{(1)}(K;\sqrt{-1})\}$
 is a basis for $V_5/{V_4}.$
 \end{thm}

We can easily see that the Vassiliev invariants $a_2(K)$,
$\sqrt{-1}F_4^{(1)}(K;\sqrt{-1})$ and $J_K^{(3)}(1)$ are additive.
If $v$ is an additive Vassiliev invariant, then, from the
coefficients of the polynomial invariants $\overline{v}$ and
$v^*$, we cannot get Vassiliev invariants other than linear
combinations of $v$ and the constant Vassiliev invariants.

Let $v$ be a Vassiliev invariant of degree $n$ and $L$ a knot.
Define $v_L^i$ to be the Vassiliev invariant defined by
$v_L^i(K)=v(L\sharp K^i)$ and define $v_L$ to be the Vassiliev
invariant defined by $v_L(K)=v(L\sharp K)$ \cite{JP2}. Then we
can see that the Vassiliev invariants obtained from the
coefficients of $\overline{v}$ and ${v}^*$ are contained in the
spans of the sets $\{v_L^i~|~L {\text{ is a knot, }} i=0, 1, 2,
\cdots, n\}$ and $\{v_L~|~L {\text{ is a knot}}\}$ respectively.

Take the trivial knot, $3_1$, $4_1$ and $5_1$ for $L$ and
$(3_1)^i$, $(4_1)^i$ and $(5_1)^i$ for $L_i$ in Theorem
\ref{vaspolth}. Then all linearly independent Vassiliev invariants
obtained by applying the $~\bar{}~$--operations and the
$^*$--operations for the non--additive Vassiliev invariants of
degree $\leq 5$ in Theorem \ref{basis} can be found as follows.
$$
\begin{aligned}
&(a_2(K))^2 \stackrel{-}{\rightarrow} \{a_2(K)\}\\
&a_4(K) \stackrel{-}{\rightarrow} \{a_2(K), (a_2(K))^2\}\\
&J_K^{(4)}(1) \stackrel{-}{\rightarrow} \{a_2(K),
(a_2(K))^2\}\\
&a_2(K)P_0^{(3)}(K;1) \stackrel{-}{\rightarrow} \{a_2(K),
J_K^{(3)}(1)\}, \quad a_2(K)P_0^{(3)}(K;1)
\stackrel{\ast}{\rightarrow} \{a_2(K)J_K^{(3)}(1)\}\\
&\quad a_2(K)J_K^{(3)}(1) \stackrel{-}{\rightarrow} \{a_2(K),
J_K^{(3)}(1)\}\\
&P_0^{(5)}(K;1) \stackrel{-}{\rightarrow}
\{a_2(K)P_0^{(3)}(K;1)\}, \quad P_0^{(5)}(K;1)
\stackrel{\ast}{\rightarrow} \{a_2(K), J_K^{(3)}(1)\}\\
&P_4^{(1)}(K;1) \stackrel{-}{\rightarrow}
\{a_2(K)P_2^{(1)}(K;1)\}, \quad P_4^{(1)}(K;1)
\stackrel{\ast}{\rightarrow} \{a_2(K), J_K^{(3)}(1)\}\\
&\quad a_2(K)P_2^{(1)}(K;1) \stackrel{-}{\rightarrow}
\{a_2(K), J_K^{(3)}(1)\}
\end{aligned}
$$
For simplicity, for each Vassiliev invariant $v,$ we unlist the
Vassiliev invariants obtained from $v^*$ if they can be obtained
from $\overline{v}$ and we also exclude the constant map $1$ whose
image is $\{1\}$ and $v$ itself in the list of Vassiliev
invariants obtained from $\overline{v}$ and $v^*$.

Thus we get the following

\begin{thm}
Let $A_n$ be a subset of the space $V_n$ of the Vassiliev
invariants of degree $\leq n$ such that
span$\mathrm{(}A_n^*\mathrm{)} = V_n$. Then $A_n$ can be chosen as follows.\\
$\mathrm{(1)}$\qua $A_0 = A_1 = \{1\},$ where $1$ denotes the constant
map with
image $\{1\}$.\\
$\mathrm{(2)}$\qua $A_2 = \{a_2(K)\}.$\\
$\mathrm{(3)}$\qua $A_3 = \{a_2(K), J_K^{(3)}(1)\}.$\\
$\mathrm{(4)}$\qua $A_4 = \{J_K^{(3)}(1), a_4(K), J_K^{(4)}(1)\}.$\\
$\mathrm{(5)}$\qua $A_5 = \{P_0^{(5)}(K;1), P_4^{(1)}(K;1),
\sqrt{-1}F_4^{(1)}(K;\sqrt{-1}), a_4(K), J_K^{(4)}(1)\}.$
\end{thm}

Let $v$ be a Vassiliev invariant of degree $n.$ In \cite{JP2}, the
authors generalized the one--variable knot polynomial invariants
$\bar v$ and $v^*$ to two--variable knot polynomial invariants
$\bar v$ and $v^*$, respectively with the same notation.

Now we want to generalize the two--variable knot polynomial
invariants $\bar v$ and $v^*$ in Theorem \ref{vaspolth}
simultaneously to a multi--variable knot polynomial invariant
$\hat v$ by unifying both $\bar v$ and $v^*$ to a multi-variable
polynomial invariant $\hat{v}$ whose proof is analogous to that of
Theorem \ref{vaspolth}. See \cite{JP2}.

Given sequences $\{L_i^{(1)}\}_{i=0}^{\infty},
\cdots,\{L_i^{(k)}\}_{i=0}^{\infty}$ of knots induced from
DD--tangles, for each knot $K$, there exists a unique polynomial
$$p_K(x_0,x_1, \cdots, x_k) \in \mathbf{F}[x_0, x_1, \cdots,
x_k]$$ such that for all $(i_0, i_1,\cdots\!, i_k) \in
\mathbb{N}^{k+1},$ $v(K^{i_0}\sharp L_{i_1}^{(1)} \sharp \cdots
\sharp L_{i_k}^{(k)}) = p_k(i_0,i_1, \cdots\!, i_k).$

Now we define a new polynomial invariant $\hat{v}\co
\{\text{knots}\}\rightarrow \mathbf{F}[x_0, \cdots, x_k]$ by
$\hat{v}(K)=p_K(x_0, \cdots, x_k)$.

Then by applying the similar argument to the case of $\bar v$ and
$v^*$ \cite{JP2}, we can see that $\hat v$ is a Vassiliev
invariant of degree $\leq n$ and the degree of each variable $x_i$
in $\hat v(K)$ is $\leq n$. Thus we get the following

\begin{thm}\label{vaspolth^}
Let $v$ be a Vassiliev invariant of degree $n$ taking values in a
numerical field $\mathbf{F}$ and let
$\{L_i^{(1)}\}_{i=0}^{\infty},
\cdots,\{L_i^{(k)}\}_{i=0}^{\infty}$ be sequences of knots induced
from DD--tangles. Then $\hat{v}\co \{\text{knots}\}\rightarrow
\mathbf{F}[x_0, \cdots, x_k]$ is a Vassiliev invariant of degree
$\leq n$ and the degree of each variable $x_i$ in $\hat v(K)$ is
$\leq n$.
\end{thm}

For a Vassiliev invariant $v$, let $C_v \co = $ \{the coefficients
of the polynomial $\hat v(K)$\}. Then, in Theorem \ref{vaspolth^},
$\hat v$ is completely determined by $C_v.$ Since the degree of
each variable in $\hat v$ is $\leq n$, we see that
$$\mathrm{span}(C_v) = \mathrm{span}(\{\hat
v(K)|_{(x_0,\cdots,x_k)=(i_0,\cdots,i_k)}~|~0\leq i_0,\cdots,i_k
\leq n\}).$$

\begin{que} Let $v$ be a Vassiliev invariant of
degree $n$. Find sequences $\{L_i^{(1)}\}_{i=0}^{\infty},$
$\cdots,$ $\{L_i^{(k)}\}_{i=0}^{\infty}$ of knots induced from
DD--tangles such that $\mathrm{span}(C_v)$ $=
\mathrm{span}(\{v\}^*)$ where $C_v$ is the set of coefficients of
the polynomial invariant $\hat v$ induced from $v$ and
$\{L_i^{(1)}\}_{i=0}^{\infty},
\cdots,\{L_i^{(k)}\}_{i=0}^{\infty}$.
\end{que}

\Addresses\recd


\begin{thebibliography}

\bibitem{BN} {\bf D Bar-Natan}, {\it On the Vassiliev knot
invariants}, Topology 34 (1995) 423--472

\bibitem{BN2} {\bf D Bar-Natan}, {\it Polynomial invariants are
polynomials}, {\tt arXiv:q-alg/9600625}

\bibitem{BL} {\bf J\,S Birman}, {\bf X-S Lin}, {\it Knot polynomials and
Vassiliev's invariants}, Invent. Math. 111 (1993) 225--270

\bibitem{JP1} {\bf M-J Jeong}, {\bf C-Y Park}, {\it Vassiliev
invariants and double dating tangle}, J. of Knot Theory and Its
Ramifications 11 (2002) 527--544

\bibitem{JP2} {\bf M-J Jeong}, {\bf C-Y Park}, {\it Vassiliev
invariants and knot polynomials}, to appear in Topology and Its
Applications

\bibitem{Ka1} {\bf T Kanenobu}, {\it An evaluation of the first derivative
of the $Q$--polynomial of a link}, Kobe J. Math. 5 (1988)
179--184

\bibitem{Ka2} {\bf T Kanenobu}, {\it Kauffman polynomials as Vassiliev link
invariants}, from: ``Proceedings of  Knots 96'', (S Suzuki, editor), World
Sci. Publ. Co. Singapore (1997) 411--431

\bibitem{Ka3} {\bf T Kanenobu}, {\bf Y Miyazawa}, {\it HOMFLY
polynomials as Vassiliev link invariants}, from: ``Knot Theory'',
Banach Center Publications 42, (V\,F\,R Jones, J
Kania--Bartoszy\'{n}ska, J\,H Przytycki, P Traczyk and V\,G Turaev,
editors), Institute of Mathematics, Polish Academy of Science,
Warszawa (1998) 165--185

\bibitem{Ka4} {\bf T Kanenobu}, {\it Vassiliev knot invariants of order $6$},
J. of Knot Theory and its Ramifications 10 (2001) 645--665

\bibitem{Kaw} {\bf A Kawauchi}, {\it A Survey of Knot Theory},
Birkh\"{a}user Verlag (1996)

\bibitem{R} {\bf D Rolfsen}, {\it Knots and Links}, Publish or
Perish Inc. (1990)

\bibitem{Sto1} {\bf A Stoimenow}, {\it Problem Session Notes},
available from his webpage: {\tt http://guests.mpim-bonn.mpg.de/alex}

\bibitem{V} {\bf V\,A Vassiliev}, {\it Cohomology of knot spaces},
           from: ``Theory of Singularities and Its Applications'',
           (V\,I Arnold, editor) Advances in Soviet Mathematics,
           Vol. 1, AMS (1990)

\bibitem{W} {\bf S Willerton}, {\it Vassiliev invariants and the Hopf
          algebra of chord diagrams}, Math. Proc. Camb. Phil. Soc.
          119 (1996) 55--65

\bibitem{W2} {\bf S. Willerton}, {\it Cabling the Vassiliev
invariants}, preprint

\end{thebibliography}
\end{document}